\documentclass[final,3p,times]{elsarticle}

\usepackage{amssymb,mathrsfs,amsmath,amsfonts,mathbbold,bbm,graphicx,float,multirow
}
\usepackage[dvipsnames,usenames]{color}

\input{undertilde}

\usepackage{natbib}
\usepackage{moreverb}
\usepackage{url}
\usepackage{grffile}
\usepackage{multirow}
\usepackage{caption}
\usepackage{graphicx}
\usepackage{float}
\usepackage{booktabs}
\usepackage{amsmath}
\usepackage{array}

\journal{Mathematica Slovaca}
\begin{document}
\captionsetup[figure]{labelfont={bf},labelformat={default},labelsep=period,name={Fig.}}


\title{On adaptivity of wavelet thresholding estimators with negatively super-additive dependent noise}

\author{Yuncai~Yu}
\author{Xinsheng~Liu\corref{cor1}}
\author{Ling~Liu}
\author{Weisi~Liu}

\cortext[cor1]{Correspondence:xsliu@nuaa.edu.cn.}

 \begin{frontmatter}

\begin{abstract}
This paper considers the nonparametric regression model with negatively super-additive dependent (NSD) noise and investigates the convergence rates of thresholding estimators. It is shown that the term-by-term thresholding estimator achieves nearly optimal and the block thresholding estimator attains optimal (or nearly optimal) convergence rates over Besov spaces. Additionally, some numerical simulations are implemented to substantiate the validity and adaptivity of the thresholding estimators with the presence of NSD noise. 
\end{abstract}
\end{frontmatter}

\section{Introduction}
\label{sec:intro}
Suppose that we observe dataset $\{Y_m\}$ by the model

\[
Y_m=g(x_m)+\varepsilon_m,~m=1,2,\cdots,n,
\tag{1.1}
\]

\par\noindent
where $x_m=m/n$ and $\varepsilon_1,\cdots,\varepsilon_n$ are identically distributed random variables defined on a probability space $(\Omega, \mathcal{F}, P)$ with zero mean and finite variance $\sigma^{2}$, $g(x)$ is an unknown function restricted to the interval $[0,1]$.

\par
There is interesting to recover $g(x)$ by wavelet method. For example, Donoho et al. (1995), Hall and Patil (1996), Hall et al. (1999), and recently Hoffmann et al. (2015) and Gao and Zhou (2016). Indeed, the adaptive estimators produced by the methods can achieve the exact minimax optimal rates, however, these adaptivity results rely on the assumption of independence noise, which is a serious restriction in model (1.1) when applied to practical problems (Wang, 1996). Some statisticians attempted to investigate the convergence rate of linear wavelet estimator in the nonparametric regression model with dependent noises (e.g. Li et al., 2008; Ding et al., 2007; Tang et al., 2018), but their rates maybe not optimal. In addition, some noises, such as long memory noise and $\rho$-mixing noise (belong to elliptically contoured family) considered by Li and Xiao (2007) and Doosti et al. (2011), respectively, involve some parameters, which are usually hard to be identified and verified. This paper considers a wide class of noise produced by NSD random sequence, whose definition based on the super-additive functions.

\par\vspace{0.2cm}\noindent
\textbf{Definition 1.1.}~~A function $\phi$: $ \mathbb {R}^{n}\rightarrow \mathbb {R}$, is called super-additive if

\[
\phi( \boldsymbol {x}\vee  \boldsymbol {y})+\phi( \boldsymbol {x}\wedge  \boldsymbol {y})\geq\phi (\boldsymbol {x})+\phi( \boldsymbol {y}),
\]

\par\noindent
for all $ \boldsymbol {x}, \boldsymbol {y}\in  \mathbb {R}^{n}$, where $``\vee"$ indicates componentwise maximum and $``\wedge"$ is for componentwise minimum.

\par\vspace{0.2cm}\noindent
\textbf{Definition 1.2.}~~A random vector $(X_1, X_2,\cdots, X_n)$ is said to be NSD if

\[
E\phi\left(X_1,X_2,\cdots,X_n\right)\leq E\phi\left(X_{1}^{*},X_{2}^{*},\cdots,X_{n}^{*}\right),
\tag{1.2}
\]

\par\noindent
where $\{X_{m}^{*}, m=1,\cdots,n\}$ are independent random variables with same marginal distribution of $\{X_{m},m=1,\cdots,n\}$ for each $m$ , and $\phi$ is a super-additive function such that the expectations in (1.2) exist.

\par\vspace{0.2cm}\noindent
\textbf{Definition 1.3.}~~A sequence of random variables $\left(X_1, X_2,\cdots, X_n\right)$ is called NSD if for all $n\geq 1$, $\left(X_1, X_2,\cdots, X_n\right)$ is NSD.

\par
The concept of NSD, which generalizes the concept of negative association (see Christofides and Vaggelatou, 2004), was proposed by Hu (2000). It is realized that many multivariate distributions possess the NSD property exhibited in practical examples, including (a) elliptically contoured distribution, (b) FGM distribution, (c) multinomial, (d) convolution of unlike multinomial, (e) multivariate hypergeometric, (f) Dirichlet, (g) Dirichlet compound multinomial, (h) negatively correlated normal distribution, (i) permutation distribution, (j) random sampling without replacement, and (k) joint distribution of ranks. Therefore, NSD has received enormous increasing attention for the potential applications in multivariate analysis and systems reliability. We refer to Hu (2000) for essential properties, Eghbal et al. (2010) for strong law of large numbers, Shen et al. (2013) for strong convergence, Wang et al. (2014) and Wang et al. (2015) for complete convergence, Shen et al. (2016) for complete moment convergence, Yu et al. (2017) for the central limit theorem. NSD samples have also been introduced to the model (1.1), and some asymptotic properties of the nonparametric regression estimators have been explored. For example, Shen et al. (2015) got the complete consistency of the weighted regression estimators by using Rosenthal-type
inequality, Wu et al. (2016) gave the convergence rate of the analogous estimators in the model (1.1) with NSD noise, Wang et al. (2018) obtained strong and weak consistency of LS estimators in the EV regression model, and Yu et al. (2019) detected the multiple change points for linear processes under NSD.

\par
The purpose of this paper is to establish the asymptotic convergence rates of the thresholding wavelet estimators in the model (1.1) with NSD noise. We demonstrate that these estimators achieve optimal and nearly optimal convergence over Besov functions class. Moreover, some simulations are implemented by \emph{R Software} to compare the block thresholding wavelet estimator with term-by-term estimator on two test functions.

\par
The remainder of this paper is organized as follows. We introduce some necessary backgrounds of thresholding estimators and state the main results in Section 2. The numerical simulations are presented to show the performances of the wavelet thresholding estimations in Section 3. Some lemmas and their proofs are provided in Section 4, and the proofs of the main theorems included in Section 5.

\section{Main results}

\subsection{Background}

Let the scaling function $\varphi$ and its associated wavelet function $\psi$ be generated from dilation equation, they are also assumed to be compactly supported and $\int\varphi=1$. In the present paper, we assume that both $\varphi$ and $\psi$ have $r$ continuous derivatives and $r$ vanishing moments, i.e., $\int_{0}^{1}x^k \varphi(x)dx=0$, $k=1,2,\cdots, r-1$, $\int_{0}^{1} x^k \psi(x) dx=0$, $k=1, 2, \cdots, r-1$.

\par
Define $\varphi_{i_0j}(x)=2^{i_0/2} \varphi(2^{i_0}x-j),~\psi_{ij}(x)=2^{i/2} \psi(2^ix-j).$ For any constant $M>0$ and a given $r$-regular wavelet $\psi$ with $r>s$, the standard Besov function space is given by

\[
B_{p,q}^s(M)=\left\{g\in B_{p,q}^s: \left(\sum\limits_{j} \left\vert\alpha_{i_0j}\right\vert^p \right)^{1/p}+
\left\{ \sum\limits_{i=i_0}^{\infty} \left[2^{i\theta} \left(\sum\limits_j \left\vert \beta_{ij}\right\vert^p \right)^{1/p} \right]^q \right\}^{1/q}\leq M, 1\leq p, q\leq \infty, 0<s<r+1\right\},
\]

\par\noindent
where $s$ is an index of regularity, $p$ and $q$ are used to specify the type of norm, $\theta=s+1/2-1/p$ and $s>1/p$.

\par
Here and below, we assume the unknown function $g(x)\in B_{p,q}^s(M) $, so that $g(x)$ can be reconstructed as

\[
g(x)=\sum\limits_{j=0}^{2^{i_0}-1} \alpha_{i_{0}j} \varphi_{i_{0}j}(x)+\sum\limits_{i\geq i_0} \sum\limits_{j=0}^{2^{i_0}-1} \beta_{ij}\psi_{ij}(x),
\tag{2.1}
\]

\par\noindent
where the coefficients are

\[
\alpha_{i_0j}=\int_0^1 g(x) \varphi_{i_0j}(x)dx,~\beta_{ij}=\int_0^1 g(x)\psi_{ij}(x)
dx.
\tag{2.2}
\]

\par
If we estimate the coefficients $\alpha_{ij}$ and $\beta_{ij}$ by $\bar{\alpha}_{ij}=n^{-1} \sum\limits_{m=1}^nY_m\varphi_{ij}(x_m)$ and $\bar{\beta}_{ij}=n^{-1} \sum\limits_{m=1}^nY_m\psi_{ij}(x_m)$, respectively. Therefore, the term-by-term thresholding estimator of (2.1) is given by

\[
\bar{g}(x)=\sum\limits_{j=0}^{2^{i_0-1}} \bar{\alpha}_{i_0j} \varphi_{i_0j}(x)+\sum\limits_{i= i_0}^{i_1} \sum\limits_{j\in\mathbb{Z}} \bar{\beta}_{ij}\psi_{ij}(x) I\left( \left\vert\bar{\beta}_{ij}\right\vert > \lambda_{0}\right),
\tag{2.3}
\]

\par\noindent
where $2^{i_0}$ is a truncation point, $i_1$ satisfies $2^{i_1-1}\leq n/\log n \leq 2^{i_1}$ and the threshold $\lambda_{0}=\sqrt{2\sigma^2\log n/n}$.

\par
Usually, term-by-term thresholding estimator produces a degree of over-smoothing. This problem can be overcome by estimating not ${\beta}_{ij}$ but its average over neighbouring coefficients (Hall and Patil (1996)). Specifically, for each resolution level $i$, we partition the integers $\{0, 1, \cdots, 2^{i}-1\}$ into consecutive, non-overlapping blocks of length $l$, that is,

\[
\Gamma_{ik}=\left\{ j: (k-1)l+1\leq j\leq kl\right\},~k\in \mathbb{Z}.
\]

\par
Let $V_i$ and $W_i$ be the spaces spanned by $\{\varphi_{ij}, j\in \mathbb{Z}\}$ and $\{ \psi_{ij}, j\in \mathbb{Z}\}$, respectively, and denote the projection operators on these spaces by $Proj_{V_i}(\cdot)$ and $Proj_{W_i}(\cdot)$. Assume the sample size $n=2^{i_2}$ and define $\hat{G}_{i_2}=n^{-1/2}\sum\limits_{m=1}^n Y_m \varphi_{i_2m}(x)$. Let the coefficients $\hat{\alpha}_{i_{0}j}$ and $\hat{\beta}_{ij}$ be given by

\[
Proj_{V_{i_0}}(\hat{G}_{i_2})=\sum\limits_{j=0}^{2^{i_0}-1}\hat{\alpha}_{i_0j}\varphi_{i_0j}~~and~~Proj_{W_{i}}(\hat{G}_{i_2})=\sum\limits_{j=0}^{2^{i}-1}\hat{\beta}_{ij}\psi_{ij}.
\]

\par
As in Hall et al. (1999, p. 42), there exist real numbers $r_{i_2m}~(m=1, 2, \cdots, n)$, such that

\[
n^{-1/2}g(m/2^{i_2}) \int\varphi=n^{-1/2}g(m/n)=\alpha_{i_2m}+r_{i_2m}.
\]

\par
Then, $\hat{G}_{i_2}(x)$ can be rewritten as

\[
\hat{G}_{i_2}=\sum\limits_{m=1}^n(\alpha_{i_2m}+r_{i_2m})\varphi_{i_2m}(x)+n^{-1/2}\sum\limits_{m=1}^{n} \varepsilon_m\varphi_{i_2m}(x).
\]

\par
Analogously, for each integer $i<i_2$, there exist real numbers $u_{ij}=\sum\limits_{m=1}^{n} r_{i_2m}<\varphi_{i_2m}, \psi_{ij}>$ and~$v_{i_0j}=\sum\limits_{m=1}^{n} r_{i_2m}<\varphi_{i_2m}, \varphi_{i_0j}>$, such that

\[
Proj_{W_i}(\hat{G}_{i_2})=\sum\limits_{j}(\beta_{ij}+u_{ij}+U_{ij})\psi_{ij}(x),
\]

\[
Proj_{V_{i_0}}(\hat{G}_{i_2})=\sum\limits_{j}(\alpha_{i_0j}+v_{i_0j}+V_{i_0j})\varphi_{i_0j}(x),
\]

\par\noindent
where

\[
U_{ij}=\dfrac{1}{\sqrt{n}}\sum\limits_{m=1}^n\varepsilon_m<\varphi_{i_2m},\psi_{ij}>,~
V_{i_0j}=\dfrac{1}{\sqrt{n}}\sum\limits_{m=1}^n\varepsilon_m<\varphi_{i_2m},\varphi_{i_0j}>.
\tag{2.4}
\]

\par
According to Li and Xiao (2010), we can obtain the block thresholding estimator

\[
\hat{g}(x)=\sum\limits_{j=0}^{2^{i_0}-1} \hat{\alpha}_{i_0j} \varphi_{i_0j}(x)+\sum\limits_{i\geq i_0}^{i_2-1} \sum\limits_{k=0}^{2^{i_0}-1} \sum\limits_{(ik)}\left(\hat{\beta}_{ij}\psi_{ij}(x)\right) I\left(\hat{B}_{ik}>\lambda^{2}\right),
\tag{2.5}
\]

\par\noindent
where $\hat{B}_{ik}=l^{-1} \sum\limits_{(ik)} \hat{\beta}_{ij}^2$ (here $\sum\limits_{(ik)}$ denotes summation over $j\in\Gamma_{ik}$), the smoothing parameter $i_0$ is chosen to satisfy $2^{i_0-1}\leq n^{1/2s+1}\leq 2^{i_0}$, the block length $l=\log n$ and the threshold $\lambda^2\geq \sigma^2 n^{-1}$.

\par
 Throughout this paper, let $C$ be a general positive constant. Put $x^+=x I(x\geq 0)$ and $x^-=-xI(x<0)$, and the inner product of $f$ and $g$ in $L^2[0, 1]$ is denoted by $<f, g>=\int fg$.

\subsection{Main theorems}

To derive our theorems, we impose a regularity condition on the noise in the model (1.1), namely, for all $k\geq 1$,

\[
v(u)=\sum\limits_{m: \vert k-m\vert\geq u} \left\vert Cov(\varepsilon_k, \varepsilon_m)\right\vert\rightarrow 0,~as~u\rightarrow \infty.
\tag{2.6}
\]

\par\vspace{0.2cm}\noindent
\textbf{Remark 2.1.}~~The condition (2.6) is easily satisfied. For example, if $v(1)<\infty$, which is the usually case, such as $\rho$-mixing sequence (required that $\sup_{m}\vert Cov(\varepsilon_m, \varepsilon_{m+u})\vert\rightarrow 0$), then $v(u)\rightarrow 0$ as $u\rightarrow \infty$. For long range dependence sequence with $\vert Cov(\varepsilon_1, \varepsilon_{1+u})\vert=Cu^{-\alpha}$, $0<\alpha\leq 1$, we have $v(u)=u^{-\alpha}$, then the condition (2.6) is satisfied as well. Particularly, independence sequence will lead to $v(1)\rightarrow 0$, which implies that the independence assumption is a serious restriction on the noise.

\par\vspace{0.2cm}\noindent
\textbf{Theorem 2.1.}~~For a given smoothing parameter $i_{0}$ and each resolution level $i$, let $\alpha_{ij}$ and $\beta_{ij}$ be given by (2.2), the noise satisfies condition (2.6), then there exists a constant $C$ such that

\[
E\left(\bar{\alpha}_{i_0j}-\alpha_{i_0j}\right)^{2} \leq C/n,
\tag{2.7}
\]

\[
E\left(\bar{\beta}_{ij}-\beta_{ij}\right)^{2} \leq C/n,
\tag{2.8}
\]

\par\noindent
and

\[
E\left(\left|\bar{\beta}_{ij}-\beta_{ij}\right|^{4}\right) \leq C2^{i}/n.
\tag{2.9}
\]

\par\vspace{0.2cm}\noindent
\textbf{Theorem 2.2.}~~In the model (1.1), assume that $\{\varepsilon_m, 1\leq m\leq n\}$ is a sequence of NSD random variables with the condition (2.6) hold. Let the wavelets $\varphi$ and $\psi$ be $r$-regular Cofflets. The term-by-term thresholding estimator $\bar{g}$ is given by (2.3) and $g(x)$ is bounded. Then for $1/p<s<r$, and $1\leq q\leq\infty$, there exists a constant $C$ such that

\[
\sup_{g\in B_{p, q}^s(M)} E\int\left(\bar{g}-g \right)^2\leq C\left(\log n/n\right)^{2s/(2s+1)}.
\]

\par\vspace{0.2cm}\noindent
\textbf{Theorem 2.3.}~~In the model (1.1), assume that $\{\varepsilon_m, 1\leq m\leq n\}$ is a bounded NSD sequence with $|\varepsilon_m|\leq L$ ($L>0$) and the condition (2.6) hold, for block thresholding estimator $\hat{g}$ given by (2.5), we have

\[
\sup_{g\in B_{p, q}^s(M)} E\int\left(\hat{g}-g \right)^2\leq\left\{ \begin{array}{ll}
Cn^{-2s/(2s+1)},&for~p\geq 2,\\
Cn^{-2s/(2s+1)}\left(\log n\right)^{(2-p)/p(1+2s)},&for~1\leq p<2.
\end{array}\right.
\]

\par\vspace{0.2cm}\noindent
\textbf{Remark 2.2.}~~Hall et al. (1999) obtain similar convergence rates over a large function space $\mathcal{H}$ with i.i.d. noise. In fact, the Theorem is still valid if we enlarge our function space $B_{p, q}^s$ by superposing the functions in it with piecewise H$\ddot{o}$lder functions similar to Hall et al. (1999), here we omit the details.

\section{Numerical study}

We take Spikes and Corner functions which were also used by Cai (1999) to illustrate the spatial adaptivity of wavelet shrinkage for independence noise as test functions. The corresponding original signals, Spikes and Corner, are assumed to be sampled at $n$ equally spaced points $x_m=m/n$, $m=1, 2, \cdots, n$. Throughout our simulations, the sample size is taken to be $n=1024$ and a low signal-to-noise ratio (SNR) is chosen (SNR=$4$). Additionally, the noise $\left\{\varepsilon_m, m=1, 2, \cdots, n\right\}$ is generated from a multivariate mixture of normal distribution with joint distribution $N (\mu_1, \mu_2, \sigma_{1}^{2}, \sigma_{2}^{2}, \rho_0)$, $\rho_0<0$, which was proven to be NSD by Yu et al. (2017). Here, $\rho_0$ is specified to depend on the signal noise level and SNR, the variances $\sigma_{1}^{2}$ and $\sigma_{2}^{2}$ are set to ensure $\rho_0<0$. These noised signals are described in Figs 1 (a) and 1 (b).

\begin{figure}[htp]
\centering
\resizebox{14cm}{66mm}{\includegraphics{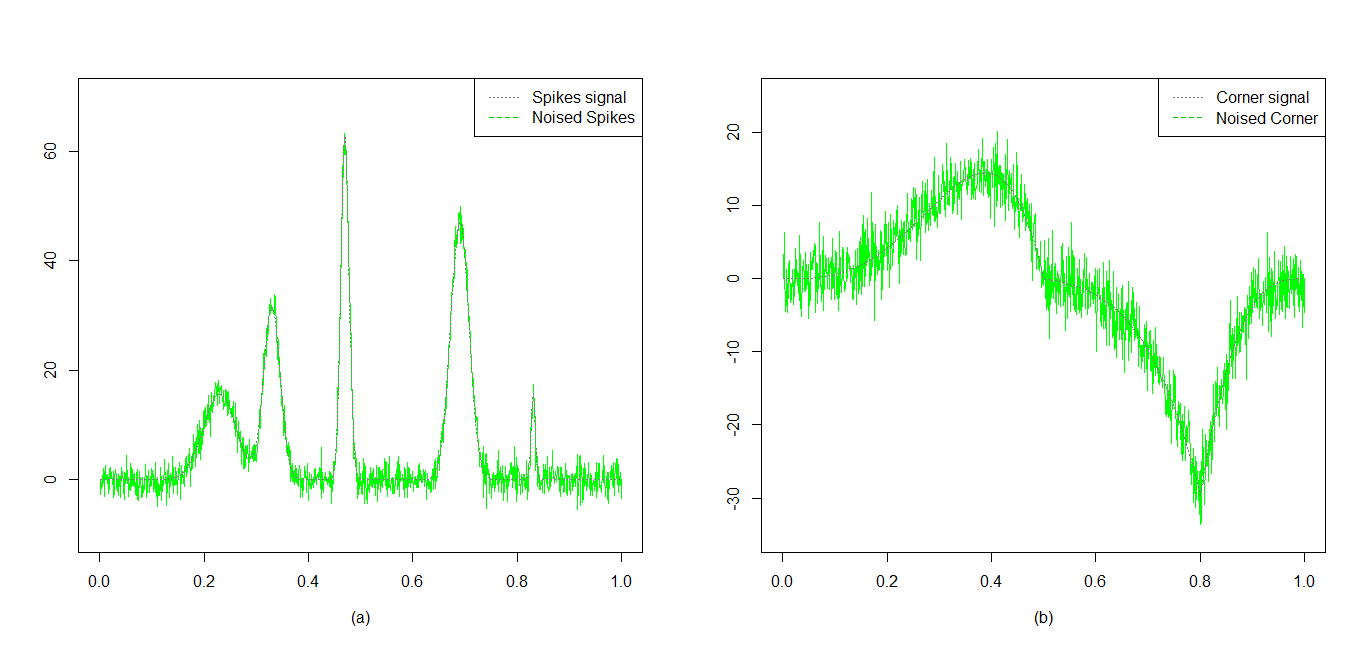}}
\caption{\label{fig:Fig1}
(a) Spikes with NSD noise; (b) Corner with NSD noise. All of the NSD noises are generated
 from a mixture of normal distribution with joint distribution $N(0, 0, 1, 9, \rho_0)$, and SNR $=4$.
}
\end{figure}

\begin{figure}[htp]
\centering
\resizebox{16cm}{110mm}{\includegraphics{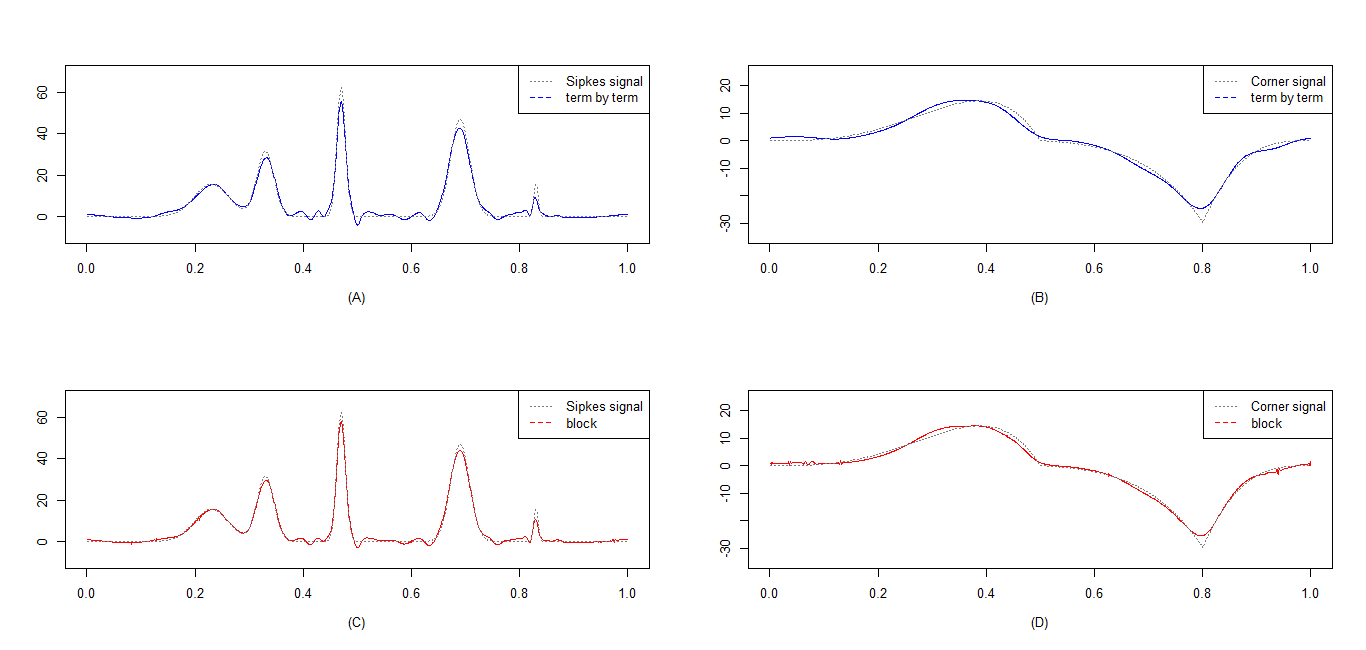}}
\caption{(A), (B) Reconstructions for Spikes and Corner using term-by-term thresholding estimator for noised Spikes and Corner with the threshold $\lambda_0=\sqrt{2\hat{\sigma}^{2}n^{-1}\log n}$; (C), (D) reconstructions for Spikes and Corner using block thresholding estimator for noised Spikes and Corner with the threshold $\lambda^2=\hat{\sigma}^{2}(x_{ik})n^{-1}$.}
\end{figure}

\par
In the case of term-by-term thresholding, we set the threshold $\lambda_0=\sqrt{2\hat{\sigma}^2n^{-1}\log n}$, where $\hat{\sigma}^2=\frac{1}{2(n-1)}\sum\limits_{m=1}^{n-1} (Y_{m+1}-Y_m)^2$. As to the block thresholding estimator, the threshold is taken to be $\lambda^2=\hat{\sigma}^2(x_{ik}) n^{-1}$, where $x_{ik}$ is a design point (one of the $x_m$s) chosen so that $2^i x_{ik}$ lies as close as possible to the middle of block $\Gamma_{ik}$, and $\hat{\sigma}^2(x_{ik})$ is an estimator of the variance of $Y_{ik}$.

\par
Figs 2 (A) and 2 (C) display that term-by-term thresholding causes some serious perturbations in the vicinity of the Spikes turning points, by contrast, block thresholding behaves relatively robust against variations of the Spikes. Figs 2 (B) and 2 (D) illustrate that block thresholding has lower bias in the vicinity of the $``$corners" (the discontinuous points of their first derivative), and provides more extensive adaptivity than term-by-term thresholding. Moreover, comparison with the term by term thresholding estimator indicates that the block thresholding estimator is superior in terms of the mean squared error(being $1.42$, $0.63$ for Spikes and Corner respectively, while $2.36$, $0.96$ for term-by-term thresholding), and thus provides extensive adaptivity to many irregularities function classes.

\section{Some lemmas}

\par\vspace{0.2cm}\noindent
\textbf{Lemma 4.1} (Hu, 2000).~If $\{X_n, n\geq 1\}$ is a NSD random sequence, we have the following properties.
\par
(a).~~Let $f_1$, $f_2$, $\cdots$ be a sequence of Borel functions all of which are non-decreasing, then $\{f_n(X_n), n\geq 1\}$ is NSD.
\par
(b).~~The sequence $\{-X_1, -X_2, \cdots, -X_n\}$ is still NSD.
\par
(c).~~If $\{X_1, X_2, \cdots, X_n\}$ is a NSD sequence, then $\{X_{\eta_1}, X_{\eta_2}, \cdots, X_{\eta_n}\}$ is NSD for any permutation $\eta=(\eta_{1},\eta_{2}, \cdots, \eta_{n})$ of $\{1,2,\cdots,n\}$.

\par\vspace{0.2cm}\noindent
\textbf{Lemma 4.2} (Yu et al., 2017). ~Suppose that $\{X_m, m\geq 1\}$ is a NSD random sequence with the condition (2.6) hold, and an array of real numbers $\{a_m, 1\leq m\leq n, n\geq 1\}$ is satisfied $\sum\limits_{m=1}^{n} a_m^2=C_0$, $C_0$ is a positive constant. Then

\[
\sigma_n^2=Var\left( \sum\limits_{m=1}^{n} a_m \varepsilon_m\right)\leq C_0 \sigma^2.
\]

\par\vspace{0.2cm}\noindent
\textbf{Lemma 4.3} (Wang et al., 2015). ~Let ${X_{n},n\leq1}$ be a NSD random sequence with mean zero and finite second moments. Denote $B_{n}^{2}=\sum\limits_{i=1}^{n}EX_{i}^{2}$. Then, for all $\epsilon > 0$, $z > 0$, and $n \geq 1$,

\[
P\left(\max\limits_{1\leq k\leq n}\left|\sum\limits_{m=1}^{k}X_{m}\right|\geq \varepsilon \right)\leq 2P\left(\max\limits_{1\leq m\leq n}\left|X_{m}\right| > z \right)
+4\exp \left\{-\frac{\varepsilon^{2}}{4(\varepsilon z+B_{n}^{2})}\right\}.
\]

\par\vspace{0.2cm}\noindent
\textbf{Lemma 4.4}. ~Let $U_{ij}$ be the random variables defined as in (2.4), take $L_{1}=\max\{\sigma,L\}$, then for all integers $i$, $k$, and real numbers $\tau\geq l\sigma^2/6n$,

\[
P\left\{ \sum\limits_{(ik)} U_{ij}^2 \geq \tau\right\}\leq C\exp\left(-\dfrac{n\tau}{8lL_{1}^2}\right).
\]

\par\vspace{0.2cm}\noindent
\textbf{Proof.} ~Let $\Delta=\left\{ w=(w_1, w_2,\cdots, w_l)\in \mathbb{R}^l : \sum\limits_{m=1}^l w_m^2=1\right\}$ be the unit sphere in $\mathbb{R}^l$. It is easy to check that for all integers $i$, $k$,

\[
\left( \sum\limits_{(ik)} U_{ij}^2\right)^{1/2} =\sup_{w\in\Delta} \sum\limits_{j=1}^l w_j U_{ij}.
\]

\par\noindent
Therefore, in order to prove Lemma 4.4, it is desired to prove that for all $t>L_{1}\sqrt{l}/\sqrt{6n}$,

\[
P\left\{ \sup_{w\in\Delta} \sum\limits_{j=1}^l w_j U_{ij}\geq t\right\}\leq C\exp\left\{-\dfrac{nt^2}{8lL_{1}^2}\right\}.
\]

\par
Consider a stochastic process $\{S_n(w), w\in \Delta\}$ as: $S_n(w)=\sum\limits_{j=1}^l w_j U_{ij}$. Denote $d_{mj}=<\varphi_{i_1m}, \psi_{ij}>$, then $S_n(w)$ can be rewritten as

\[
S_n(w)=\dfrac{1}{\sqrt{n}} \sum\limits_{m=1}^{n} \varepsilon_m \sum\limits_{j=1}^l w_j d_{mj}.
\]

\par
Define $S_n^+(w)=\dfrac{1}{\sqrt{n}} \sum\limits_{m=1}^{n} \sum\limits_{j=1}^l \left( w_j d_{mj} \right)^{+} \varepsilon_m$, $S_n^-(w)=\dfrac{1}{\sqrt{n}} \sum\limits_{m=1}^{n} \sum\limits_{j=1}^l \left( w_j d_{mj} \right)^{-} \varepsilon_m$. For $\delta>0$, we will consider the cases $0<|\delta S_{n}^{+}(w)|\leq 1$ and $1<|\delta S_{n}^{+}(w)|$ in the following proof.

\par
\textbf{Case 1:} For $0<|\delta S_{n}^{+}(w)|\leq 1$. Since $E\varepsilon_{m}=0$, we have

\[
\begin{array}{lll}
E\sup\limits_{w\in\Delta}  \exp\left\{ \delta S_n^+(w) \right\}
&\leq & 1+\sum\limits_{k=2}^{\infty}\frac{E\sup\limits_{ w\in\Delta}(\delta S_n^+(w))^{k}}{k!} \\
&\leq & 1+\delta^{2}E\sup\limits_{ w\in\Delta}(S_n^+(w))^{2}\left\{ \frac{1}{2!}+\frac{1}{3!}+ \cdots \right\} \\
&\leq & 1+\delta^{2}E\sup\limits_{ w\in\Delta}(S_n^+(w))^{2}\\
&\leq & \exp\left\{\delta^{2}E\sup\limits_{ w\in\Delta}(S_n^+(w))^{2}\right\}.
\tag{4.1}
\end{array}
\]

\par
In the view of $\sum\limits_{m=1}^l w_m^2=1$, by Schwarz's inequality, it follows

\[
\begin{array}{lll}
\exp\left \{ \delta^2 E\sup\limits_{w\in\Delta} \vert S_n^+(w)\vert ^2\right\}
&\leq& \sup \limits_{w\in \Delta}\prod\limits_{m=1}^{n}
\left( \exp\left\{E\delta\dfrac{1}{n}  \sum\limits_{j=1}^l \vert w_j d_{mj} \vert^2 \varepsilon_m^+ \right\}\right) \\
&\leq& \exp\left\{\sigma ^2 \delta^2 \cdot \dfrac{1}{n}  \left(\sum\limits_{m=1}^n \vert d_{mj} \vert^2 \right)^{1/2}\right\}.
\end{array}
\]

\par
Note that $d_{mj}^2=1$, if $m=j$; and $d_{mj}^2=0$, otherwise. For a fixed $j~\left(j=1,2, \cdots,l\right)$, there only exists a counterpart $m=j$ such that $d_{mj}^2=1$ (otherwise $d_{mj}^2=0$), which implies that $\sum\limits_{m=1}^n \left( \sum\limits_{j=1}^l \vert d_{mj}\vert^2 \right)^{1/2}=\sqrt{l}$, and this leads to

\[
\begin{array}{lll}
\exp\left\{ \delta^2 E\sup\limits_{w\in\Delta} \left\vert S_n^+(w)\right\vert ^2\right\} \leq \exp\left\{ \sigma^2\delta^2\sqrt{l}/n\right\}.
\tag{4.2}
\end{array}
\]

\par
\textbf{Case 2:} For $|\delta S_{n}^{+}(w)|> 1$. According to the properties (a) and (b) in Lemma 4.1, for every $|\varepsilon_{m}|\leq L$, we obtain by taking $\phi(x)=e^x$ in (1.2) that

\[
\begin{array}{lll}
E\sup\limits_{w\in\Delta}  \exp\left\{ \delta S_n^+(w) \right\}
&\leq& E\exp\left \{ \delta^2\sup\limits_{w\in\Delta} \vert S_n^+(w)\vert ^2\right\} \\
&\leq& \sup \limits_{w\in \Delta}\prod\limits_{m=1}^{n}
E\left( \exp\left\{\delta^2 \dfrac{1}{n}  \sum\limits_{j=1}^l \vert w_j d_{mj} \vert^2 \varepsilon_m^2 \right\}\right) \\
&\leq& \sup \limits_{w\in \Delta}\prod\limits_{m=1}^{n}
\left( \exp\left\{\delta^2 L^2 \dfrac{1}{n}  \sum\limits_{j=1}^l \vert w_j d_{mj} \vert^2 \right\}\right) \\
&\leq& \exp\left\{\delta^2 L^2 \cdot \dfrac{1}{n}  \left(\sum\limits_{m=1}^n \vert d_{mj} \vert^2 \right)^{1/2}\right\}\\
&=& \exp\left\{ \delta^2 L^2\sqrt{l}/n\right\}.
\tag{4.3}
\end{array}
\]

\par
Thus take (4.1), (4.2) and (4.3) together to give

\[
E\sup\limits_{w\in\Delta}  \exp\left\{ \delta S_n^+(w)\right\} \leq \exp\left\{ L_{1}^2\delta^2\sqrt{l}/n\right\}.
\]

\par
By Markov's inequality, for every $t>0$ and $\delta>0$,

\[
P\left\{\sup\limits_{w\in\Delta} S_n^+(w)\geq t\right\}\leq \dfrac{ E \exp\left\{ \sup\limits_{w\in\Delta}  S_n^+(w)\right\}  } {\exp\left\{\delta t\right\}}\leq \exp\left\{  -\delta t+L_{1}^2\delta^2\sqrt{l}/n\right\}.
\]

\par\noindent
Hence, for $\delta=\dfrac{nt}{2L_{1}^2 \sqrt{l}}$, we get that

\[
P\left\{\sup\limits_{w\in\Delta} S_n^+(w)\geq t\right\}\leq C\exp\left\{ -t^2n/4L_{1}^2\sqrt{l} \right\}.
\]

\par
Analogously, for $S_n^-(w)>0$, we have

\[
P\left\{\sup\limits_{w\in\Delta} S_n^-(w)\geq t\right\}\leq C\exp\left\{ -t^2n/4L_{1}^2\sqrt{l} \right\}.
\]

\par\noindent
Note that $\sup\limits_{w\in\Delta} S_n(w)\leq \sup\limits_{w\in\Delta} S_n^+(w)+\sup\limits_{w\in\Delta} S_n^-(w)$, then

\[
\begin{array}{lll}
P\left\{\sup\limits_{w\in\Delta} S_n(w)\geq t\right\} &\leq& P\left\{\sup\limits_{w\in\Delta} S_n^+(w)\geq t/2\right\} + P\left\{\sup\limits_{w\in\Delta} S_n^-(w)\geq t/2\right\}\\
&\leq& C\exp\left\{ -t^2n/8L_{1}^2\sqrt{l} \right\}.
\end{array}
\]

\par\noindent
This completes the proof of Lemma 4.4.

\section{Proof of the theorems}

\par\vspace{0.2cm}\noindent
\textbf{Proof of Theorem 2.1.}~~It is easy to verify that $\bar{\alpha}_{i_0j}$ is an unbiased estimator of $\alpha_{i_0j}$ by

\[
E\bar{\alpha}_{i_0j}=E\left(n^{-1} \sum\limits_{m=1}^nY_m\varphi_{i_0j}(x_m)\right)=\int\varphi_{i_0j}(x)g(x)dx=\alpha_{i_0j},
\]

Note that $\sum\limits_{m=1}^n<\varphi_{i_{0}m}, \varphi_{i_0j}>^2=1$, based on Lemma 4.2, it yields

\[
E\left(\bar{\alpha}_{i_0j}-\alpha_{i_0j}\right)^{2}
= Var\left(n^{-1} \sum\limits_{m=1}^nY_m\varphi_{i_0j}(x_m)\right)= n^{-1} Var\left(n^{-1/2} \sum\limits_{m=1}^n\varepsilon_m\varphi_{i_0j}(x_m)\right)
\leq C_{0}\sigma^{2}/n = C/n.
\]

Obviously, (2.8) hold since these inequalities still hold for replacing $\bar{\alpha}_{ij}$, $\alpha_{ij}$ by $\bar{\beta}_{ij}$ and $\beta_{ij}=\int_{0}^{1}g(x)\psi_{ij}(x)dx$, respectively. Recall that $g(x)$ is bounded, using Cauchy-Schwarz's inequality, we obtain

\[
|\beta_{ij}| \leq \int_{0}^{1}g(x)|\psi_{ij}|dx \leq C\int_{0}^{1}|\psi_{ij}|dx \leq C\left(\int_{0}^{1}|\psi_{ij}|dx\right)^{2}=C.
\tag{5.1}
\]

\par\noindent
According to $E\varepsilon_{m}=0$, it follows

\[
|\bar{\beta}_{ij}| \leq \sup_{1\leq m \leq n}\psi_{ij}(x)EY_{m} \leq \sup_{1\leq m \leq n}\psi_{ij}(x_{m})|g(x_{m})| \leq C\sup_{1\leq m \leq n}
\psi_{ij}(x_{m}) \leq C2^{i/2}.
\tag{5.2}
\]

Combining (5.1) and (5.2), we derive that

\[
E\left|\bar{\beta}_{ij}-\beta_{ij}\right|^{4} \leq \left(\bar{\beta}_{ij}+\beta_{ij}\right)^{2}E\left(\bar{\beta}_{ij}-\beta_{ij}\right)^{2} \leq C\left(1+2^{i/2}\right)^{2}n^{-1}=C2^{i}/n.
\]

\par\vspace{0.2cm}\noindent
\textbf{Proof of Theorem 2.2.}~~We will divide the proof of Theorem 2.2 into several parts. By the orthogonality of $\varphi$ and $\psi$, the quadratic risk can be decomposed as

\[
E\Vert \bar{g}-g\Vert_2^2=R_1+R_2+R_3,
\]

\par\noindent
where
\par
$R_1=\sum\limits_{j=0}^{2^{i_0}-1} E\left(\bar{\alpha}_{i_0j}-\alpha_{i_0j}\right)^2,$
\par
$R_2=\sum\limits_{i=i_1}^{\infty} \sum\limits_{j=0}^{2^i-1} \beta_{ij}^2,$
\par
$R_{3}=\sum\limits_{i=i_0}^{i_1-1}\sum\limits_{j=0}^{2^{i}-1}E\left(\bar{\beta}_{ij}I\left(\bar{\beta}_{ij}>\lambda_{0}\right)-\beta_{ij}\right)^{2}$.

\par\noindent
The reminder of the proof consists of bounding $R_1$, $R_2$, $R_3$.
\par
Bound for $R_1$: From $2s/(2s+1)<1$ and (2.7) in Theorem 2.1, we get

 \[
R_1\leq C\sum\limits_{j=0}^{2^{i_0}-1}1/n \leq C2^{i_0}/n \leq C\left(\log n/n\right)^{2s/(2s+1)}.
\tag{5.3}
\]

\par
Bound for $R_2$: For $p\geq2$, we have $B_{p,q}^{s}(M) \subseteq B_{2,\infty}^{s}(M)$, then

\[
R_{2}^{\dag}\leq C\sum\limits_{i=i_1}^{\infty}2^{-2is}\leq C2^{-i_1s} \leq C\left(\log n/n\right)^{2s/(2s+1)}.
\tag{5.4}
\]

\par
For $1\leq p \leq 2$, we have $B_{p,q}^{s}(M) \subseteq B_{2,\infty}^{s+1/2-1/p}(M)$, by $(s+1/2-1/p) > s/(2s+1)$, one can see that

\[
R_{2}^{*}\leq C\sum\limits_{i=i_1}^{2^{i_0}-1}2^{-2i(s+1/2-1/p)}\leq C2^{i_1(s+1/2-1/p)} \leq C\left(\log n/n\right)^{2s/(2s+1)}.
\tag{5.5}
\]

Combining (5.4) and (5.5), we obtain

\[
R_{2}\leq C\left(\log n/n\right)^{2s/(2s+1)}.
\tag{5.6}
\]

 \par
 Bound for $R_3$: Note that $\left\{\left|\bar{\beta}_{ij}\right|<\lambda_{0}, \left|\beta_{ij}\right| \geq 2\lambda_{0}\right\} \subseteq \left\{\left|\bar{\beta}_{ij}-\beta_{ij}\right|>\lambda_{0}/2\right\},
 \left\{\left|\bar{\beta}_{ij}\right|<\lambda_{0}, \left|\beta_{ij}\right|\geq 2\lambda_{0}\right\}\subseteq \left\{\left|\beta_{ij}\right|\leq
 2\left|\bar{\beta}_{ij}-\beta_{ij}\right|\right\}$ and $\left\{\left|\bar{\beta}_{ij}\right|\geq \lambda_{0},
  \left|\beta_{ij}\right|< \lambda_{0}/2\right\}\subseteq \left\{\left|\bar{\beta}_{ij}-\beta_{ij}\right|>\lambda_{0}/2\right\}$. We can write $R_{3}$ as

 \[
 \begin{array}{lll}
 R_{3}&=&\sum\limits_{i=i_0}^{i_1-1}\sum\limits_{j=0}^{2^{i}-1}
 E\left\{\left(\bar{\beta}_{ij}-\beta \right)^{2}
 I\left(\bar{\beta}_{ij}\geq\lambda_{0}\right)\right\}+
 \sum\limits_{i=i_0}^{i_1-1}\sum\limits_{j=0}^{2^{i}-1}E\left\{\beta_{ij}^{2}I\left(\bar{\beta}_{ij}<\lambda_{0}\right)\right\}\\
 &\leq& 2\sum\limits_{i=i_0}^{i_1-1}
 \sum\limits_{j=0}^{2^{i}-1}E\left\{\left(\bar{\beta}_{ij}-\beta \right)^{2}I\left(\bar{\beta}_{ij}\geq\lambda_{0}\right)\right\}
 +\sum\limits_{i=i_0}^{i_1-1}\sum\limits_{j=0}^{2^{i}-1}E\left\{\left(\bar{\beta}_{ij}-
 \beta \right)^{2}I\left(\bar{\beta}_{ij}\geq\lambda_{0} \right)
 I\left(\beta_{ij}\geq\lambda_{0}/2\right)\right\} \\
 &&+\sum\limits_{i=i_0}^{i_1-1}\sum\limits_{j=0}^{2^{i}-1}E\left\{\beta_{ij}^{2}I\left(\|\beta_{ij}\|< \lambda_{0}\right)I\left(\|\beta_{ij}\|< 2\lambda_{0}\right)\right\}\\
 &=:& 2R_{31}+R_{32}+R_{33}.
 \tag{5.7}
 \end{array}
 \]

\par
Firstly, we bound $R_{31}$ by Cauchy-Schwarz's inequality

\[
R_{31}\leq \sum\limits_{i=i_0}^{i_1-1}\sum\limits_{j}\left(E\left(\bar{\beta}_{ij}
-\beta_{ij}\right)^{4}\right)^{1/2}\left(P\left(|\bar{\beta}_{ij}
-\beta_{ij}|>\lambda_{0}/2\right)\right)^{1/2}.
\]

Next, we consider the bound of the probability $P(|\bar{\beta}_{ij}
-\beta_{ij}|>\lambda_{0}/2)$. Note that

\[
\left|\bar{\beta}_{ij}
-\beta_{ij}\right|=\left|n^{-1}\sum\limits_{m=1}^{n}Y_{m}\psi_{ij}(x_m)-n^{-1}\sum\limits_{m=1}^{n}g(x_{m})\psi_{ij}(x_m)\right|
=n^{-1}\left|\sum\limits_{m=1}^{n}\psi_{ij}(x_{m})\varepsilon_{m}\right|.
\]

\par

Without loss of generality, we assume $\psi_{ij}(x_{m})>0$. From the property (c) in Lemma 4.1, all positive $\varepsilon_m$ can be moved to the first $k$ terms, hence

\[
 n^{-1}\left|\sum\limits_{m=1}^{n}\psi_{ij}(x_{m})\varepsilon_{m}\right|
\leq \max\limits_{1\leq m \leq k}n^{-1}\left|\sum\limits_{m=1}^{k}\psi_{ij}(x_{m})\varepsilon_{m}\right|
= \max\limits_{1\leq m \leq k}n^{-1}\left|\sum\limits_{m=1}^{k}\psi_{ij}(x_{m})\varepsilon_{m}^{+}\right|
\leq \max\limits_{1\leq m \leq k} Cn^{-1}2^{i}\left|\sum\limits_{m=1}^{k}\varepsilon_{m}\right|.
\]

\par\noindent
Since $|\bar{\beta}_{ij}-\beta_{ij}|\leq |\bar{\beta}_{ij}|+|\beta_{ij}|\leq C2^{i/2}$, $B_{n}^{2}=E\sum\limits_{m=1}^{n}\varepsilon_{m}^{2}=n\sigma^{2}$ and $2^{i} \leq n$, apply Lemma 4.3, with $z=\sqrt{2\sigma}2^{i/2}$, to obtain

\[
\begin{array}{lll}
P\left(\left|\bar{\beta}_{ij}
-\beta_{ij}\right|>\lambda_{0}/2\right)&\leq &P\left(\max\limits_{1\leq m\leq k}\left|\sum\limits_{m=1}^{k}\varepsilon_{m}\right|> n\lambda_{0} 2^{-i-1}\right) \\
&\leq & 4\exp\left\{-\frac{n^{2}\lambda_{0}^{2}2^{-(i+1)}}{4\left(n
\lambda_{0} \sqrt{2\sigma} 2^{-(i+1)/2}+n\sigma^{2}\right)}\right\} \\
& \leq& 4\exp\left\{-\frac{\sqrt{\log n} \cdot 2^{-(i+1)/2}}{4C\sqrt{n}}\right\}\\
&\leq & C\exp \left\{-\frac{\sqrt{\log n}}{4}\right\} = O\left(n^{-1/4}\right).
\end{array}
\]

\par\noindent
From (2.9) in Theorem 2.1, it follows

\[
R_{31}\leq C \sum\limits_{i=i_0}^{i_1-1}\sum\limits_{j}\left(2^{i}/n\right)^{1/2}
\left(1/n^{-4}\right)^{1/2}\leq C n^{-2}\sum\limits_{i=i_0}^{i_1-1}2^{i}\leq C n^{-2}2^{i_1}
 \leq C\log n/n \leq C\left(\log n/n\right)^{2s/(2s+1)}.
 \tag{5.8}
\]

In order to bound $R_{32}$, we choose $i_{3}$ to satisfy $2^{i_3-1}\leq \left(n/\log n\right)^{1/(2s+1)} \leq 2^{i_3}$, using (2.8) in Theorem 2.1, we get

 \[
 \begin{array}{lll}
 R_{32} &\leq& C/n \cdot \sum\limits_{i=i_0}^{i_1-1}\sum\limits_{j}
  I\left(|\beta_{ij}|>\lambda_{0}/2\right) \\
 &\leq& C/n \cdot \sum\limits_{i=i_0}^{i_{3}-1}\sum\limits_{j}I\left(
 |\beta_{ij}|>\lambda_{0}/2\right)+C/n \cdot \sum\limits_{i=i_{3}}^{i_1-1}\sum\limits_{j}I\left(
 |\beta_{ij}|>\lambda_{0}/2\right) =: R_{321}+R_{322}.
 \end{array}
 \]

For the term $R_{321}$, one can show that

\[
R_{321}\leq C/n \cdot \sum\limits_{i=i_0}^{i_{3}-1}2^{i}
\leq C2^{i_{3}-1}/n \leq C\left(\log n/n\right)^{2s/(2s+1)}.
\]

To bound the term $R_{322}$, we consider the cases $p\geq 2$ and $1\leq p<2$ separately. For $p\geq2$, by Markov's inequality, it follows

\[
R_{322}^{\dag}\leq C\lambda_{0}^{-2}/n\sum\limits_{i=i_{3}}^{\infty}\sum
\limits_{j}\beta_{ij}^{2} \leq C/n \cdot 2^{i_{3}}\leq C\left(\log n/n\right)^{2s/(2s+1)}.
\]

\par
For $1\leq p \leq 2$, note that $I(|\beta_{ij}|>\lambda_{0}/2)\leq
C|\beta_{ij}|^{p}/{\lambda_{0}}^{p}$ and $(2s+1)(2-p)/2+(s+1/2-1/p)p = 2s$, then

\[
R_{322}^{*}\leq C \left(n\lambda_{0}^{p}\right)^{-1}\sum\limits_{i=i_{3}}^{i_1}\sum\limits_{j}\left|\beta_{ij}\right|^{p}\leq Cn^{-(2-p)/2}\sum\limits_{i=i_{3}+1}^{\infty}2^{-i(s+1/2-1/p)}\leq C\left(\log n/n\right)^{2s/(2s+1)}.
\]

Hence, $R_{32}$ can be bounded by

\[
R_{32}\leq C\left(\log n/n\right)^{2s/(2s+1)}.
\tag{5.9}
\]

\par
In the view of (5.7), (5.8) and (5.9), it is sufficient to bound the term $R_{33}$, which may be written as

\[
\begin{array}{lll}
R_{33}&\leq& \sum\limits_{i=i_0}^{i_1-1}\sum\limits_{j}
\beta_{ij}^{2}I\left(|\beta_{ij}|<2\lambda_{0}\right) \\
&=& \sum\limits_{i=i_0}^{i_{3}-1}\sum\limits_{j}
\beta_{ij}^{2}I\left(|\beta_{ij}|<2\lambda_{0}\right)
+\sum\limits_{i=i_{3}}^{i_1-1}\sum\limits_{j}
\beta_{ij}^{2}I\left(|\beta_{ij}|<2\lambda_{0}\right):=R_{331}+R_{332}.
\end{array}
\]

\par
Based on $\beta_{ij}^{2}I(|\beta_{ij}|<2\lambda_{0})\leq 4\lambda_{0}^{2}$, we get that

\[
R_{331}\leq C\sum\limits_{i=i_0}^{i_{3}}2^{i}\lambda_{0}^{2}\leq C2^{i_{3}}\log n/n\leq C\left(\log n/n\right)^{2s/(2s+1)}.
\]

By performing the same operation for the term $R_{332}$ as $R_{322}$, we can show that $R_{332}\leq C\left(\log n/n\right)^{2s/(2s+1)}$. Then

\[
R_{3}\leq C\left(\log n/n\right)^{2s/(2s+1)}.
\tag{5.10}
\]

\par
Thus the Theorem 2.2 follows from (5.3), (5.6) and (5.10).

\par\vspace{0.2cm}\noindent
\textbf{Proof of Theorem 2.3.}~~Similarly to the proof of Theorem 2.2, we decompose the quadratic risk as several parts

\[
E\Vert \hat{g}-g\Vert_2^2=T_1+T_2+T_3+T_4,
\]

\par\noindent
where

\par
$T_1=\sum\limits_{j=0}^{2^{i_0}-1} E\left(\hat{\alpha}_{i_0j}-\alpha_{i_0j}\right)^2,$
\par
$T_2=\sum\limits_{i=i_{2}}^{\infty} \sum\limits_{j=0}^{2^i-1} \beta_{ij}^2,$
\par
$T_3=\sum\limits_{i=i_0}^{i_{2}-1}\sum\limits_{k\in \mathbb{Z} } P\left(\hat{B}_{ik}\leq\lambda^{2}\right) \sum\limits_{(ik)}  \beta_{ij}^2,$
\par
$T_4=\sum\limits_{i=i_0}^{i_{2}-1}\sum\limits_{k\in \mathbb{Z} } E\left\{ I(\hat{B}_{ik}>\lambda^{2}) \sum\limits_{(ik)} \left(\hat{\beta}_{ij}-\beta_{ij}\right) ^2\right\}.$
\par
In the view of $\sum\limits_{m=1}^n<\varphi_{i_{2}m}, \varphi_{ij}>^2=1$ and $\sum\limits_{m=1}^n<\varphi_{i_{2}m}, \psi_{i_0 j}>^2=1$, by Lemma 4.2, it yields

\[
EU_{ij}^2\leq \sigma^2/n, ~EV_{i_0j}^2\leq \sigma^2/n.
\tag{5.11}
\]

\par\noindent
According to the definition of $\hat{\alpha}_{i_0 j}$ and (5.11), one gets

\[
T_1=E\left\Vert Proj_{V_{i_0}}\left(\hat{G}_{i_{2}}-g\right)\right\Vert_2^2=\sum\limits_{j=0}^{2^{i_0}-1} v_{i_0 j}^2+ \sum\limits_{j=0}^{2^{i_0}-1}E V_{i_0 j}^2\leq \sum\limits_{m=1}^{n}r_{i_{2}m}^{2}+ 2^{i_0}\sigma^2/n\leq Cn^{-2s/(1+2s)}.
\]

\par
Furthermore, by (5.6), we can show that

\[
T_2=\sum\limits_{i=i_{2}}^{\infty} \sum\limits_{j=0}^{2^i-1} \beta_{ij}^2 \leq C n^{-2s/(1+2s)}.
\]

\par
From Lemma 4.7 in Li et al. (2010, p.1119), we see that

\[
T_3\leq \left\{ \begin{array}{ll}
C n^{-2s/(1+2s)},&~p\geq 2,\\
C\left(\log n\right)^{(2-p)/p(1+2s)} n^{-2s/(1+2s)},&~1\leq p<2.
\end{array}\right.
\]

\par
To complete the proof of Theorem 2.3, we focus on Bounding $T_4$ by appropriate rates. Firstly, we write $T_4$ as

\[
\begin{array}{lll}
T_4&=&\sum\limits_{i=i_0}^{i_{2}-1}\sum\limits_{k} E\left\{ I\left(\hat{B}_{ik}>\lambda^{2}\right) \sum\limits_{(ik)} (U_{ij}+u_{ij}) ^2\right\}\\
&\leq& 2\sum\limits_{i=i_0}^{i_{2}-1}\sum\limits_{k} E\left\{ I\left(\hat{B}_{ik}>\lambda^{2}\right) \sum\limits_{(ik)} u_{ij}^2\right\}+2\sum\limits_{i=i_0}^{i_{2}-1}\sum\limits_{k} E\left\{ I\left(\hat{B}_{ik}>\lambda^{2}\right) \sum\limits_{(ik)} U_{ij}^2\right\}\\
&=:& 2T_4^{'}+2T_4^{''}.
\end{array}
\tag{5.12}
\]

\par
Then from (5.11), we see that $T_4^{'}$ is bounded by

\[
T_4^{'}\leq \sum\limits_{i=i_0}^{i_{2}-1} \sum\limits_{k} \sum_{(ik)} u_{ij}^2 \leq  \sum\limits_{i=i_0}^{i_{2}-1} \sum\limits_{j=0}^{2^i-1} u_{ij}^2  \leq \sum\limits_{m=1}^n r_{i_{2} m}^2\leq Cn^{-2s/(1+2s)}.
\tag{5.13}
\]

\par
Finally, we consider the bound of the term $T_{4}^{''}$. For $p\geq 2$, we choose $i_{4}$ to satisfy $2^{i_{4}-1}\leq n^{1/2s+1}\leq 2^{i_{4}}$, and denote $B_{ik}=l^{-1} \sum\limits_{(ik)} (\beta_{ij}+u_{ij})^2$. Thus $T_4^{''}$ can be divided into three parts:

\[
\begin{array}{lll}
T_4^{''}&=&\sum\limits_{i=i_0}^{i_{4}}\sum\limits_{k} E\left\{ I\left(\hat{B}_{ik}>\lambda^{2}\right) \sum\limits_{(ik)} U_{ij}^2\right\} +\sum\limits_{i=i_{4}+1}^{i_{2}-1}\sum\limits_{k} E\left\{ I\left(\hat{B}_{ik}>\lambda^{2}\right) I\left(B_{ik}>\lambda^{2}/2\right) \sum\limits_{(ik)} U_{ij}^2\right\}\\
&&+\sum\limits_{i=i_{4}+1}^{i_{2}-1}\sum\limits_{k} E\left\{ I\left(\hat{B}_{ik}>\lambda^{2}\right) I\left(B_{ik}\leq \lambda^{2}/2\right) \sum\limits_{(ik)} U_{ij}^2\right\}\\
&=:& T_{41}+T_{42}+T_{43}.
\end{array}
\]

\par\noindent
From (5.11), we have

\[
T_{41}\leq \sum\limits_{i=i_0}^{i_{4}} \sum\limits_{k} \sum\limits_{(ik)} E U_{ij}^2 = \sum\limits_{i=i_0}^{i_{4}} \sum\limits_{j=0}^{2^i-1} \sigma^2 n^{-1} \leq \sigma^2 n^{-1}  \sum\limits_{i=i_0}^{i_{4}} 2^i \leq Cn^{-2s/(1+2s)}.
\tag{5.14}
\]

\par\noindent
Recall that $\lambda^{2}\geq \sigma^2 n^{-1}$, thus $T_{42}$ is bounded by (5.11) that

\[
T_{42}\leq 2\sum\limits_{i=i_{4}+1}^{i_{2}-1} \sum\limits_{k}B_{ik} \lambda^{-2}  \sum\limits_{(ik)} E U_{ij}^2 =\dfrac{1}{24} \sum\limits_{i=i_{4}+1}^{i_{2}-1} \sum\limits_{j=0}^{2^i-1} \left(\beta_{ik}+u_{ij}\right)^2 \leq \sum\limits_{i=i_{4}+1}^{i_{2}-1} 2^{-2si} + \sum\limits_{m=1}^{n} r_{i_{2} m}^2 \leq Cn^{-2s/(1+2s)}.
\tag{5.15}
\]

\par\noindent
To bound the term $T_{43}$, we appeal to the Lemma 5.1 in Hall et al. (1999, p. 45), which implies that

\[
\left\{\hat{B}_{ik}>\lambda^{2}\right\} \cap \left\{B_{ik}\leq \frac{\lambda^{2}}{2}\right\} \subseteq \left\{ \sum\limits_{(ik)} U_{ij}^2 \geq \dfrac{l}{12} \lambda^{2}\right\} =\left\{ \sum\limits_{(ik)} U_{ij}^2 \geq \dfrac{\sigma^2 l}{6n} \right\}.
\]

\par\noindent
Hence

\[
T_{43}\leq  \sum\limits_{i=i_{4}+1}^{i_{2}-1} \sum\limits_{k=0}^{2^i-1} \int_{\sigma^2 l/6n}^{\infty} P\left\{\sum\limits_{(ik)} U_{ij} ^2\geq \tau \right\}d\tau.
\]

\par\noindent
Form Lemma 4.4, it follows that

\[
T_{43}=C  \sum\limits_{i=i_{4}+1}^{i_{2}-1} \sum\limits_{k=0}^{2^i-1} \int_{\sigma^2 l/6n}^{\infty} \exp\left\{ -\dfrac{n\tau}{8\sqrt{l}L_{1}^2}\right\}d \tau= \dfrac{C}{n} \sum\limits_{i=i_{4}+1}^{i_{2}-1} 2^i \exp\left\{-\sqrt {l}/48\right\}=C\exp\left\{-\sqrt{l}/48\right\}.
\]

\par
Since $l=\log n$, we have $T_{43}\leq O\left(n^{-\eta}\right)$ for all $\eta>0$. Putting this result and (5.12), (5.13), (5.14), (5.15) together, one can see that $T_4\leq Cn^{-2s/(2s+1)}$ when $p\geq 2$.

\par
For the case $1\leq p<2$. We treat $T_4$ similarly as (4.16) in Li et al. (2010, p. 1123), then

\[
T_{4} \leq C\left(\log n\right)^{(2-p)/p(1+2s)} n^{-2s/(1+2s)}.
\]

\par
Consequently, $T_4$ is bounded by

\[
T_4\leq \left\{ \begin{array}{ll}
C n^{-2s/(1+2s)},&~p\geq 2,\\
C\left(\log n\right)^{(2-p)/p(1+2s)} n^{-2s/(1+2s)},&~1\leq p<2.
\end{array}\right.
\]

\par
This completes the proof of Theorem 2.3.

\section*{Acknowledgements}
This paper is supported by Natural Science Foundation of China (No. 61374183); the Postgraduate Research and Practice Innovation Program of Jiangsu Province (No. KYCX$19\_0149$).\\

\vspace{0.5cm}
{
\makeatletter
\renewcommand\@biblabel[1]{}
\renewenvironment{thebibliography}[1]
     {\section*{\refname}%
      \@mkboth{\MakeUppercase\refname}{\MakeUppercase\refname}%
      \list{\@biblabel{\@arabic\c@enumiv}}%
           {\settowidth\labelwidth{\@biblabel{#1}}%
            \leftmargin\labelwidth
            \advance\leftmargin\labelsep
            \itemindent-\labelsep
            \@openbib@code
            \usecounter{enumiv}%
            \let\p@enumiv\@empty
            \renewcommand\theenumiv{\@arabic\c@enumiv}}%
      \sloppy
      \clubpenalty4000
      \@clubpenalty \clubpenalty
      \widowpenalty4000%
      \sfcode`\.\@m}
     {\def\@noitemerr
       {\@latex@warning{Empty `thebibliography' environment}}%
      \endlist}
\makeatother

\makeatletter
\renewenvironment{thebibliography}[1]
{\section*{\refname}%
\@mkboth{\MakeUppercase\refname}{\MakeUppercase\refname}%
\list{\@biblabel{\@arabic\c@enumiv}}%
{\settowidth\labelwidth{\@biblabel{#1}}%
\leftmargin\labelwidth \advance\leftmargin\labelsep
\advance\leftmargin by 2em%
\itemindent -2em%
\@openbib@code
\usecounter{enumiv}%
\let\p@enumiv\@empty
\renewcommand\theenumiv{\@arabic\c@enumiv}}%
\sloppy \clubpenalty4000 \@clubpenalty \clubpenalty
\widowpenalty4000%
\sfcode`\.\@m} {\def\@noitemerr
{\@latex@warning{Empty `thebibliography' environment}}%
\endlist}
\makeatother

\small


\nocite{*}
\bibliographystyle{amsplain}

}

Yuncai Yu, State Key Laboratory of Mechanics and Control of Mechanical Structures, Institute of Nano Science and Department of Mathematics, Nanjing University of Aeronautics and Astronautics, Nanjing 210016, China, E-mail address: yuyuncai@nuaa.edu.cn\\

Xinsheng Liu$^{*}$, State Key Laboratory of Mechanics and Control of Mechanical Structures, Institute of Nano Science and Department of Mathematics, Nanjing University of Aeronautics and Astronautics, Nanjing 210016, China, E-mail address: xsliu@nuaa.edu.cn\\

Ling Liu, Department of Information Science and Technology, Donghua University, Shanghai 201600, China, E-mail address: lgliu0@163.com\\

Weisi Liu, State Key Laboratory of Mechanics and Control of Mechanical Structures, Institute of Nano Science and Department of Mathematics, Nanjing University of Aeronautics and Astronautics, Nanjing 210016, China, E-mail address: panda.si@hotmail.com

\end{document}